\input amstex
\loadbold
\documentstyle{amsppt}
\vsize=7.5 in
\hsize=6.0 in
\parskip\smallskipamount
\overfullrule=0pt
\magnification\magstep1
\pageno=1

\document
\baselineskip=18pt    \define\eh{e^{i\theta}}
\define\th{\theta}  \define\cc{\Bbb C}

\define\bd{\partial}  \define\be{\overline{\bd}} 
    \define\zb{\overline{z}}

\define\hf{{1\over 2}} \define\ze{\zeta}

\def\icc{\int_{\cc}} \def\rr{\Bbb R}               
\define\loc{\text{loc}}  \define\1op{{1\over \pi}}

\topmatter
\title Some conjectures about integral means of $\partial f$ and 
$\overline{\partial} f.$ 
\endtitle \author Albert Baernstein II and Stephen J. Montgomery-Smith \endauthor
\address Washington University, St. Louis, Missouri 63130,
\newline University of Missouri, Columbia, Missouri 65211-0001 \endaddress
\email al\@math.wustl.edu, Stephen\@math.missouri.edu \endemail
\thanks The first author was supported in part by NSF grants DMS-9206319 and 
DMS-9501293. The second author was supported in part by NSF grant 
DMS-9424396 \endthanks \endtopmatter

\noindent 1. \bf The problems. \rm  In this note we shall discuss some 
conjectural integral inequalities which are related to quasiconformal 
mappings, singular integrals, martingales and the calculus of variations. For 
a function $f:\cc  
\rightarrow \cc,$ denote the formal complex derivatives by
$$ \bd f = {{\bd f} \over {\bd z}} = \hf ({{\bd f} \over {\bd x}} -
i{{\bd f} \over {\bd y}}), \quad \be f = {{\bd f} \over {\be z}} =   
\hf ({{\bd f} \over {\bd x}} +  i{{\bd f} \over {\bd y}}).$$

\noindent Define a function $L:\cc \times \cc \rightarrow \rr$ by
$$\align L(z, w) & = |z|^2 - |w|^2,\quad \text{if} \quad |z| + |w| \leq 1, \\
   & = 2|z| - 1, \quad \text{if} \quad |z| + |w| > 1.\endalign$$

\proclaim {Conjecture 1} (V. \v Sver\'ak) Let $f\in \dot W^{1,2}(\cc, \cc).$ Then $$\icc L(\bd f, 
\be f)  
\geq 0. \tag1.1$$  \endproclaim

Here we denote by $\dot W^{1,2}(\cc, \cc)$ denotes the ``homogeneous" Sobolev 
space of complex valued locally integrable functions   
in the plane whose distributional first derivatives are in 
$L^2$ on the plane. Integrals without specified variables are understood to be with 
respect to Lebesgue measure. 

Since $\bd (\overline{f}) = \overline{\be f},$  Conjecture 1 is true if and 
only if (1.1) always holds when $L(\bd f,\be f)$ in the integral is replaced by
$L(\be f, \bd f).$ 
                                            
The function $L,$ with a plus 1 added to the right hand side, was introduced 
by 
Burkholder [Bu4], [Bu5, p.20]. In his setting, the variables $z$ and $w$ are 
taken from an arbitrary Hilbert space. It appears independently in work of 
\v Sver\'ak [Sv1], who considered the question, as yet unresolved, of whether 
functions belonging to a certain class which contains a function 
naturally associated to $L$ are ``quasiconvex". As will  
be explained in \S 5, quasiconvexity of this function implies (1.1). 

For $p\in (1, \infty),$ set 
$$p^* = \max\; (p, p'), \quad \text{where}\quad {1\over p} + {1\over {p'}} = 
1.$$ Again following Burkholder [Bu5, p.16], [Bu3, p.77], [Bu4, p.8], 
define functions ${\Phi}_p :\cc \times \cc \rightarrow \Bbb R\,$ by
$$  {\Phi}_p(z,w) = {\alpha}_p ((p^* - 1)|z| - |w|)\,(|z| + |w|)^{p-1},\quad 
{\alpha_p} = p(1 - {1\over {p^*}})^{p-1}.$$

For $1<p<2$ and $z, w\in \cc$ one calculates $$\int_0^{\infty} t^{p-
1}L({z\over t}, {w\over t})\,dt =  {\beta}_p\,{\Phi}_p(z,w),
\quad {\beta}_p = (\hf\,  p\,(2-p)\,{\alpha}_p)^{-1}.\tag1.2a$$           

Set $M(z,w) = L(z,w) - (|z|^2 - |w|^2) = (|w|^2 - (|z|-1)^2) 1_{(|z|+|w| > 
1)}.$ Then, for $2<p<\infty,$ 
$$\int_0^{\infty} t^{p-1}M({w\over t}, {z\over t})\,dt =  
{\gamma}_p\,{\Phi}_p(z,w),\quad {\gamma}_p = (\hf\, p\,(p-1)\,(p-2)\,{\alpha}_p )^{- 
1}.\tag1.2b$$  

For $f\in \dot W^{1,2}(\cc),$  one sees by
Fourier transforms or otherwise that $$\int_{\cc} |\bd f|^2 - |\be f|^2 = 0.$$

If Conjecture 1 is true, the last three identities imply the truth of

\proclaim {Conjecture 2} (R. Ba\~nuelos - G.Wang) For 
$f\in \dot W^{1,p}(\cc, \cc)$ holds
$$\icc {\Phi}_p(\bd f, \be f) \geq 0, \quad 1<p<\infty.$$ \endproclaim

Ba\~nuelos and Wang arrived at Conjecture 2 in the course of their work [BW1].
The conjecture is stated as Question 1' in [BL, \S 5], where the reader can 
find other questions and comments related to the present paper.

From [Bu3, p.77] follows the inequality
$${\Phi}_p (z,w) \leq  (p^* - 1)^p |z|^p - |w|^p ,\quad z,w \in \cc, \quad 1<p<\infty.$$

Thus, if Conjecture 2 is true, then so is the following conjecture, which is 
due to T.Iwaniec [I1, I3]. 

\proclaim {Conjecture 3} (T.Iwaniec) For $f\in \dot W^{1,p}(\cc, \cc)$ holds
$$\icc |\be f|^p \leq  (p^* - 1)^p \icc|\bd f|^p, \quad 1<p<\infty.$$ \endproclaim

Like Conjecture 1, each of Conjectures 2 and 3 is true if and only the same 
inequality always holds when $\bd f$ and $\be f$ are interchanged in the 
corresponding integral. 

All three conjectured inequalities are sharp, if they are true. Let $f(z) = 
cz$ for $|z|<1,\quad f(z) = c/\zb\;$ for $|z|>1,$ where $c$ is a nonzero 
complex constant. Simple calculations show that 
equality holds for $f$ in Conjecture 1. Equality holds for $f$ in Conjecture 2 
when $1< p \leq 2\;$ and for $\overline {f}$ when $2\leq p <\infty.$  
In section 6, we'll see that equality holds in 
Conjectures 1 and 2 for a large class of functions in $\dot W^{1,2}$ and 
$\dot W^{1,p},$ respectively. 

By contrast, it seems plausible that when $p\ne 2$ equality never holds in 
Conjecture 3. To construct sequences which saturate the upper bound in 
Conjecture 3, take $p\in (1, \infty)$ and $\alpha \in (0,1/p).$  Define
$f_{\alpha}(z) = z|z|^{-2\alpha},\; \text{if}\;|z|\leq 1,\;
f_{\alpha}(z) = 1/{\zb},\;\text{if}\; |z|\geq 1.$ One computes that $\icc |\bd 
f_{\alpha}|^p / \icc |\be f_{\alpha}|^p  
\rightarrow (p -1)^p$ as $\alpha \rightarrow 1/p.$  Thus, 
$\icc |\bd (\overline{f_{\alpha}})|^p / \icc |\be\,(\overline{f_{\alpha}})|^p 
\rightarrow (p -1)^{-p}$  
as $\alpha \rightarrow 1/p.$ From these two relations, and the 
interchangeability of $\bd f$, $\be f$ in the conjectures, it follows that the 
constant on the right hand side of Conjecture 3 must be at least 
$\max^p (p-1, {1\over {p-1}}) = (p^*-1)^p.$ 

Here are three reasons we find these conjectures of interest.

(a) Truth of Conjecture 3 would imply in the limiting case 
$p\rightarrow\infty$ a notable recent theorem of Astala [As1] about area 
distortion of quasiconformal mappings in the plane.   

(b)  Let $S$ be the singular integral operator in the plane defined by
$$Sf(z) = -{1\over {\pi}} \icc {f(\ze)\over {(z-\ze)^2}}\,|d\ze|^2.\tag1.3$$

Truth of Conjecture 3 would show that the norm of $S$ on $L^p(\cc, \cc)$ is 
precisely $p^* - 1.$ 

(c) Falsity of Conjecture 1 would prove, for $2\times 2$ matrix valued 
functions, a 
conjecture of Morrey in the calculus of variations which asserts that 
rank one functions are not necessarily quasiconvex. \medskip

In Sections 2-5 we'll elaborate on statements (a), (b), and (c). In 
Sections 6-8  we'll present some evidence in favor of the conjectures.

We are grateful to Professors Astala, Ba\~nuelos, and Iwaniec for helpful 
communications, especially to Professor Iwaniec for sharing some of his 
unpublished notes with us. Thanks go also to N. Arcozzi, D.Burkholder, R. 
Laugesen and the referee for corrections and comments on the first version of 
the manuscript. The first author thanks the organizers of the 
Uppsala conference for their marvelous hospitality and efficiency. Above all, 
he thanks Matts and Agneta for many years of inspiration and friendship. 

\bigskip

\noindent 2. \bf Area distortion by quasiconformal mappings. \rm The integral in 
(1.3) is a Cauchy principal value. The operator $S$ is sometimes called the 
Beurling-Ahlfors transform. The general theory of such operators, as developed 
by Calder\'on, Zygmund, and others, is presented, 
for example, in [S]. Among its consequences are the facts that $S$ is 
bounded on $L^p$ for  
$1<p<\infty,$ and that $S$ may be defined  via Fourier multipliers by 
$$(Sf)\sphat \,(\xi) = (\overline{\xi}/\xi)\hat f(\xi),\quad \xi \in \cc.\tag2.1$$  

Thus, $S$ acts isometrically on $L^2(\cc, \cc)$. It follows also from (2.1) 
that for appropriate functions $f$ we have $$S(\be f) = \bd f.\tag2.2$$ 

Because of (2.2), the Ahlfors-Beurling operator plays an important role in the 
theory of quasiconformal mapping in the plane. See, for example, [LV].
A homeomorphism $F:\cc 
\rightarrow \cc$ is said to be $K$ quasiconformal, $K\geq 1,$ if $F\in 
W_{\loc}^{1,2}(\cc,\cc),$ and if $|\be F(z)| \leq k|\bd F(z)|$ for a.e $z\in 
\cc,$ where $k = (K-1)/(K+1).$ In the 1950's, Bojarski [Bo 1,2] applied  
the recently proved $L^p$-boundedness of $S$ to prove that partial 
derivatives of $K-qc$  
maps, which a priori belong to $L^2_{\loc},$ belong in fact to $L^p_{\loc}$ for 
some $p>2$ which depends only on $K.$ Via H\"older's inequality, this 
enhanced integrability leads to 
an inequality for the distortion of area by qc maps. One way to state the area 
distortion property is as follows: If $F(0) = 0$ and $F(1) = 1,$ then for 
all measurable sets $E\subset (|z|<1),$  
$$ |F(E)| \leq C |E|^{\kappa},\tag2.3$$
where $|\cdot|$ denotes Lebesgue measure, and $C$ and $\kappa$ depend only on 
$K.$  

Gehring and Reich [GR] conjectured in 1966 that the best possible, i.e. 
smallest, $\kappa$ for which (2.3) is valid should be $\kappa = 1/K.$ 
Prototypical conjectured extremals were the radial stretch maps  
$$F_K(z) = z|z|^{{1\over K}-1}.$$ 
$F_K$ is $K-qc$ for each $ K\geq1,$ and satisfies $ |F_K(B)| 
=  {\pi}^{1-{1\over K}} |B|^{1/K}$ for balls $B$ centered at the origin.  

The Gehring-Reich conjecture withstood many assaults before it was 
finally proved in the 1990's by Astala [As1], by means of 
very innovative considerations involving holomorphic dynamics 
and thermodynamical formalism. Eremenko and Hamilton [EH] gave a shorter proof 
of the conjecture using a distillation of Astala's ideas. More background and 
related results can be found in the survey [As2].  In [N] and [AsM], the 
distortion results are applied to problems about ``homogenization" of composite 
materials. 

To continue our story requires a backup. The weak 1-1 and $L^2$ boundedness of 
$S$ imply existence of absolute constants $c$ and $\alpha$ such that for 
all $E\subset (|z| < 1),$ 
$$\int_{|z|<1}|S(1_E)| \leq c |E|\,\log({{\alpha}\over {|E|}}). $$

Gehring and Reich showed that their area distortion conjecture is more or less 
equivalent to  
proving that the smallest $c$ for which some $\alpha$ exists is $c = 1.$ Let 
$||S||_p$ denote the norm of $S$  
acting as an operator from $L^p(\cc, \cc)$ into itself. Iwaniec [I1]
found that ``c=1" is implied by 
$$\liminf_{p\rightarrow \infty}{1\over p}||S||_p  = 1.$$

This implication, together with the examples $f_{\alpha}$ in Section 1,
led Iwaniec to Conjecture 3, which, as noted in (b) at the 
end of section 1, can be restated as  
 $$||S||_p = (p^*-1),\quad 1<p<\infty. \tag2.4$$  

Thus, if Conjecture 3 is true, it could be regarded as a significantly stronger  
form of Astala's area distortion theorem.  \bigskip

\noindent 3. \bf Norms of singular integral operators and martingale 
transforms.  \rm The 
prototypical singular integral operator is the Hilbert transform $H$, defined 
for functions on $\rr$ by $$  Hf(x) = \1op \int_{\rr} {{f(y)}\over {x-
y}}\,dy.$$ 

When the Fourier transform $\hat f$ is defined as in [S], we have the 
multiplier equation  $(Hf)\sphat \,(\xi) =  
i{{\xi}\over {|\xi |}} \hat f(x).$ Thus $H$ is an isometry on $L^2(\rr).$ 
M.Riesz, in 1927, proved that $H$ is in fact bounded on $L^p(\rr),$ for 
$1<p<\infty.$ The sharp $L^p$ bounds for real valued functions were found by 
Pichorides [Pi] in 1972. Let $||H||_p$ denote the norm of $H$ acting on 
$L^p(\rr, \rr).$ Recall that $p^* = \max (p, p').$

\proclaim {Pichorides's Theorem} $\quad ||H||_p = \cot\,{{\pi}\over 
{2p^*}},\quad 1<p<\infty.$ \endproclaim                            

Let $F$ be the harmonic extension of $f+i Hf$ to the upper half plane. Then 
$F$ is holomorphic. If 
$\phi$ is subharmonic on the range of $F$ then $\phi\circ F$ is subharmonic. 
Pichorides proved his theorem by making a good choice for $\phi.$ According 
to [G], the best constant for Riesz's theorem was found independently by 
B.Cole, whose work established a generalized version of the theorem in 
the context of ``Jensen measures" on uniform algebras. 

Verbitsky [V], and a little later Ess\'en [E], gave a shorter proof of 
Pichorides's theorem  
by finding an even better choice for the subharmonic function $\phi.$ Grafakos
[Gr] found a still shorter proof. Verbitsky and Ess\'en also 
proved sharp bounds for the analytic projection operator $I+iH,$ where $I$ 
denotes the identity, acting as an operator from $L^p(\rr, \rr)$ to $L^p(\rr, 
\cc)$:  $$||I + iH||_p = \csc {{\pi}\over {2p^*}}.\tag3.1$$

As explained in [Pe], the norm of $H$ acting on $L^p(\rr, \cc)$ is still 
$\cot\,{{\pi}\over {2p^*}}.$ But for $I+iH$ acting on $L^p(\rr, \cc)$ the norm 
is apparently not known. See [Pe], [KV]. The sharp weak 1-1 constant for 
$H$ acting on $L^1(\rr, \rr)$ was found by B.Davis in 1974, but seems to be not
known for $H$ acting on $L^1(\rr, \cc).$ See [Pe], [Bu4, p.6] for discussion.  

The Beurling-Ahlfors operator $S$ furnishes one analogue of the Hilbert 
transform in dimension 2. But the most basic generalizations to higher 
dimensions of the Hilbert transform  are the \it Riesz transforms \rm  
$R_j,\;j=1,...,n.$  In terms of Fourier multipliers, they are 
defined by $$(R_jf)\sphat \,(\xi) = i{{{\xi}_j}\over {|\xi |}} \hat 
f(\xi),\quad \xi  
\in {\rr}^n,$$ and in terms of integrals by convolution with the kernel 
$C_n x_j/|x|^{n+1},$ where \newline $C_n = \Gamma({{n+1}\over 2}){\pi}^{- 
{{n+1}\over 2}} .$ 

Let $||R_j||_p$ denote the norm of $R_j$ acting on $L^p({\rr}^n, \rr)$ or
$L^p(\rr, \cc)$. 
T.Iwaniec and G.Martin [IM 3], proved that $R_j$ has the same 
norm as $H$.

\proclaim {Iwaniec-Martin Theorem} $\quad||R_j||_p = \cot\,{{\pi}\over 
{2p^*}},\quad 1<p<\infty.$ \endproclaim  

Iwaniec and Martin use the method of rotations to show that $||R_j||_p$  
is bounded above by Pichorides's constant. The lower estimate for $||R_j||_p$ is 
proved by a simple but clever ``transference'' argument.  

Ba\~nuelos and Wang [BW1] obtained another proof of $\quad||R_j||_p \leq
\cot\,{{\pi}\over {2p^*}}\,$ 
as a consequence of a `` $\cot\,{{\pi}\over {2p^*}}$ theorem" they proved for 
transformations of certain stochastic integrals. In addition, they  proved 
that (3.1) holds when $H$ is replaced by one of the $R_j.$ Apparently, the 
$L^p$ norms of operators such as $I \bigoplus R_1 \bigoplus R_2: L^p({\rr}^n, 
\rr)\rightarrow  
L^p({\rr}^n, {\rr}^3)$ remain unknown when $n\geq 2$ and $p\ne 2.$ 

Arcozzi [Ar1], see also [Ar2], [ArL], carried the martingale methods over to compact 
manifolds, Lie groups, and Gauss space. Among other things, he proves that for 
suitable definitions of Riesz  
transforms $R$ on n-spheres and on certain compact Lie groups again hold 
$||R||_p = \cot\,{{\pi}\over {2p^*}}$ and $||I + iR||_p = \csc
\,{{\pi}\over {2p^*}}.$ For general compact Lie groups, Arcozzi proves that one 
has the upper bounds $||R||_p \leq \cot\,{{\pi}\over {2p^*}}$ and $||I + 
iR||_p \leq \csc \,{{\pi}\over {2p^*}}.$

Recall that Conjecture 3 can be stated as $||S||_p 
= p^* - 1$, where $S$ is the Ahlfors-Beurling operator (1.3). This conjecture 
differs from the established result $||R_j||_p = \cot\,{{\pi}\over 
{2p^*}}$ in two respects: (i) The kernel $z^{-2}$ for $S$ is even, so the 
method of rotations is not applicable. (ii) The kernel for $S$ is complex-
valued. 

At least two sharp inequalities for $S$, not involving $L^p$, do 
exist. See [EH] and [I2]. Additional evidence for Conjecture 3 is provided in [AIS], where it is 
shown that the operator $I-S\mu$ is invertible in $L^p$ for all 
functions $\mu\in L^{\infty}(\cc,\cc)$ with $||\mu||_{L^{\infty}}\leq k$ if and 
only if $k<1/p^*.$ 

In [IM 1,2,3], Iwaniec and Martin introduce operators $S_n$ which operate on 
functions $f:{\rr}^n\rightarrow \Lambda,\,$ where $\Lambda$ is the usual 
Grassmann algebra of ${\rr}^n.$ The operator $S_2$ can be identified with $S$.
Iwaniec and Martin conjecture that in all dimensions, one still has $||S||_p =
p^*-1.$ They point out that such a result would have strong consequences for 
the regularity theory of quasiregular maps in ${\rr}^n.$ For subsequent work 
on $S_n,$ see [BL]. The survey [I4] discusses sundry related subjects in 
n dimensions.\bigskip

\noindent 4. \bf Differential subordination. \rm The Ba\~nuelos-Wang work 
continues a line of sharp constant investigation initiated 
by Burkholder in the late 1970's. The survey [Bu5] contains a good  
bibliography for this rich and varied body of work. Here we'll confine 
discussion to the parts most pertinent to Conjectures 1 and  2.  

Following [Bu5, p.16], let $\{f_n\}$ and $\{g_n\},\; n\geq 0,$ be Hilbert 
space valued martingales with respect to the same filtration on some 
probability space $(\Omega, \Cal F, P).$ Denote the corresponding difference 
sequences by $\{d_n\}$ and $\{e_n\},$ so  
that $$f_n = \sum_{k=0}^n d_k,\quad g_n = \sum_{k=0}^n e_k.$$

For $p\geq 1,$ write $||f_n||_p$ for the $L^p$ norm of $f_n$ with respect to 
$P.$ Then $||f||_p \equiv \lim_{n\rightarrow \infty} ||f_n||_p $ exists.

\proclaim {Burkholder's Theorem}  Suppose that, for all $k\geq 0$ and $P-
a.e \;\omega \in \Omega,$
$$|e_k(\omega)| \leq |d_k(\omega)|. \tag4.1$$
Then $$||g||_p \leq (p^* - 1) ||f||_p,\tag4.2$$ and the constant $p^* - 1$ is best 
possible. If $0<||f||_p < \infty,$ then equality occurs in (4.2) if and only 
if $p=2$ and equality holds a.e in (4.1) for all $k\geq 0.$ \endproclaim

Burkholder's proof of (4.2) is similar in spirit to the proofs of 
Pichorides's theorem: He shows that, with ${\Phi}_p$ the function we defined in 
section 1, the sequence of expectations $E \Phi_p(f_n,g_n)$ is nondecreasing 
for $n\geq 0,$ with
$E \Phi_p(f_0 ,g_0) \geq 0.$ It follows that the analogue of 
Conjecture 2 holds in Burkholder's setting, and hence so does the analogue 
(4.2) of Conjecture 3.

By 1984, Burkholder [Bu1] had proved that (4.1) implies (4.2) whenever the martingales 
are real valued. Among other features of the proof in [Bu1] is a reduction to 
the case when $e_k = {\epsilon}_k d_k,$  where the ${\epsilon}_k$ are 
constants, each of which is $1$ or $-1$. Then $f$ and $g$ are somewhat like 
conjugate harmonic functions, and the transform $f\rightarrow g$ can be 
viewed as an analogue of the Hilbert transform. Extension of $(4.1) \implies 
(4.2)$ to the Hilbert space valued case followed in 1988. 

When martingales $f$ and $g$ as above satisfy (4.1), Burkholder says that 
$g$ is differentially subordinate to $f.$ He introduced also, in [Bu4], 
the notion of differentially subordinate harmonic functions. If $u$ and $v$ 
are Hilbert  
space valued harmonic functions on a domain $D\subset {\rr}^n,$ then $v$ is 
said to be differentially subordinate to $u$ if $|\nabla v(x)| \leq |\nabla 
u(x)|$ at each $x\in D.$ For example, when $n=2$ each member of a pair of 
conjugate harmonic functions is differentially subordinate to the other. 
It turns out that ${\Phi_p}(u,v)$ is 
subharmonic in $D$ when $v$ is differentially subordinate to $u.$ Let $x_0\in 
D,$ and let $\mu$ denote the harmonic measure of $D$ at $x_0.$ If we assume 
also that $|v(x_0)|\leq |u(x_0)|,$ then the subharmonicity leads to the 
inequality $$ ||v||_p \leq (p^* - 1)||u||_p, \quad 1<p<\infty,\tag4.3$$
where the $L^p$ norm is taken with respect to $\mu.$

It is not known if $p^* - 1$ is best possible in (4.3). Pichorides's theorem 
implies that the best constant $c_p$ must satisfy $c_p\geq \cot\,{{\pi}\over 
{2p^*}}.$ Related papers about differential subordination include [Bu6], 
[C1], and [C2].  

The function $L$ apparently first appears in [Bu4]. A related function 
appears in [Bu2, (8)]. In [Bu4], Burkholder proves  
integral inequalities analogous to our Conjecture 1 for differentially 
subordinate martingales 
$f,\,g$ and differentially subordinate harmonic functions $u,\,v$ with 
$|v(x_0)|\leq |u(x_0)|.$  These inequalities, valid in the Hilbert space valued 
case, imply, for $n\geq 0,$ the weak-type inequalities
$$P(|f_n| + |g_n|\geq 1) \leq 2||f||_1,\quad \mu(|u| + |v| \geq 1)\leq 
2||u||_1.\tag4.4$$

The constant $2$ in (4.4) is best possible even when $|f_n| + |g_n|$ is 
replaced by $|g_n|\,$ and $|u| + |v|\,$ by $|v|.\,$ See [Bu4, p.11] and  [Bu6, 
Remark 13.1].

Burkholder proved versions of his discrete parameter martingale 
results for some continuous parameter martingales. Ba\~nuelos and Wang [BW1,2] 
and Wang [W] extended the theory to cover a wider class of
continuous parameter martingales. Theorem 1 of [BW1], about real-valued 
differentially subordinate martingales, when combined with a probabilistic 
representation of the Riesz transforms due to Gundy and Varopoulos [GV], 
leads to the new proof of the upper bound in the Iwaniec-Martin theorem and 
the proof of (3.1)  
with $H$ replaced by $R_j$ mentioned at the end of the section 3.  
Theorem 2 of [BW1] leads to the following $L^p$ estimates for the Beurling-
Ahlfors transform $S.$

\proclaim {Ba\~nuelos-Wang Theorem} $$||S||_p \leq 4(p^* - 1),\quad 
1<p<\infty.\tag4.5$$ \endproclaim 

The constant 4 in (4.4) is the smallest known at present which works for all 
$p.$ Recall that  $||S||_2 = 1,$ and that Conjecture 3 may be stated as 
$||S||_p = p^* - 1.$ 

Here are some ideas from the proof of (4.5). Take a 
rapidly decreasing smooth function $f:\cc \rightarrow \cc.$ Extend $f$ 
to a harmonic function in the  half space ${\rr}_+^3,$ denoted also by $f.$  
Then, using It\^o's formula,
$$f(B_0) = \int_{-\infty}^0 \nabla f(B_s)\cdot dB_s,$$
where $B_s$ is the ${\rr}_+^3$ - valued Gundy-Varopoulos ``background 
radiation process",   and $\nabla f \cdot dB_s$ is the complex number obtained by 
splitting $\nabla f$ into real and imaginary parts, then taking dot products 
in ${\rr}^3.$ For suitable functions $A$ whose values are complex $3\times 3-$ 
matrices, define random variables $A*f$ by 
$$A*f = \int_{-\infty}^0 (A(B_s) \nabla f(B_s))\cdot dB_s.$$ 

When the limit $0$ in the integrals is replaced by $t\in (-\infty, 0),$ one 
obtains complex valued martingales to which the extended Burkholder theory is 
applicable. Let $||A|| = \sup \{|A(z) v|:
z\in {\rr}_+^3,\, {v\in \cc}^3,\, |v|\leq 1\},$ where $|.|$ denotes the Euclidean 
norm in ${\cc}^3.$ From Theorem 2 of [BW1] follows 
$$||A*f||_p \leq ||A|| (p^* - 1) ||f(B_0)||_p.\tag4.6$$   

Now $S$ can be expressed in terms of Riesz transforms: $S={R_2}^2 - {R_1}^2 +
2i R_1 R_2.$ If $A$ is taken to be the constant matrix 
$$ A = \pmatrix 0 & 0 & 0 \\ 0 & 2 & 2i \\ 0 & 2i & -2 \endpmatrix, $$
it turns out that $$Sf(z) = E(A*f\,| \,B_0 = z),\quad z\in \cc. \tag4.7$$

The conditional expectation operator in (4.7) is a contraction on $L^p$ when 
$p\geq 1,$ and the distribution of $B_0$ on $\cc$ is Lebesgue measure. These 
facts together with (4.6) yield                                                             
$$ ||Sf||_p \leq ||A*f||_p \leq ||A||\;(p^* - 1) ||f(B_0)||_p = (p^* - 
1)||A||\;||f||_p.\tag4.8$$ 

Calculation gives $||A|| = 4.$ So (4.5) follows from (4.8).  

How is the pair $\bd f,\,\be f$ like a pair of 
differentially subordinate martingales or harmonic functions? 
That, it seems, is what we really need to know to get the full conjectured 
result $ ||Sf||_p \leq (p^* - 1)||f||_p$ by the route of this section. \bigskip

\noindent 5. \bf Quasiconvex and rank one convex functions. \rm Let 
${\rr}^{nm}$ denote the set of all $m\times n$ matrices with real 
coefficients. A function  
$\Psi :{\rr}^{nm} \rightarrow \rr$ is said to be rank one convex on 
${\rr}^{nm}$
if for each $A,\;B \in {\rr}^{nm}$ with rank $B = 1$ the function
$$h(t) \equiv \Psi(A+tB),\quad t\in \rr,$$ is convex. $\Psi$ is said to 
be quasiconvex on ${\rr}^{nm}$ if it is locally integrable and for each $A\in 
{\rr}^{nm},$ each bounded domain $D\subset {\rr}^n$  
and each compactly supported Lipschitz function $f:D\rightarrow {\rr}^m$ holds
$${1\over {|D|}}\int_D \Psi(A + \nabla f) \geq \Psi(A).$$

If $n=1$ or $m=1$ then $\Psi$ is quasiconvex or rank one convex if and only if
it is convex. If $m\geq 2$ and $n\geq 2$, then convexity $\implies$ 
quasiconvexity $\implies$ rank one convexity. See [D1], where one finds also 
a discussion of polyconvexity, a property which lies in between 
convexity and quasiconvexity. Additional relevant works include [DDGR],
[AD], [D2], [Sv1], [Sv2], and [Sv3].

Morrey [M, p.26] conjectured in 1952 that rank one convexity does not imply 
quasiconvexity  
when $m$ and $n$ are both $\geq 2$. \v Sver\' ak [Sv2], in 1992, proved 
that Morrey's  
conjecture is correct if $m\geq 3\,$ and $n\geq 2.$ The cases $m=2,\,n \geq 2$ 
remain open.  

Define $\alpha : {\rr}^{2,2} \rightarrow {\cc}^2$ by 
$\alpha (\left [\smallmatrix a&b\\ c&d \endsmallmatrix\right]) = (z_1, 
z_2),$ where $$z_1 = \hf((a+d) + i(c-b)),\quad z_2 = \hf((a-d) +i(c+b)).$$
For $f: \cc \rightarrow \cc,$ represent $\nabla f$ as a real $2\times 2$ 
matrix in the  
usual way:  $\nabla f = \left [\smallmatrix u_x & u_y \\ v_x & v_y 
\endsmallmatrix\right],$ where $u$ and $v$ are the real 
and imaginary parts of $f.$ Then $\alpha (\nabla f) = (\bd f, \be f).$

Recall that the Burkholder-\v Sver\'ak  function $L:{\cc}^2 \rightarrow \rr$ is 
defined by 
$$\align L(z_1,z_2) &= |z_1|^2 - |z_2|^2,\quad \text{if}\quad |z_1| + |z_2| \leq 
1,\\ & = 2|z_1| - 1, \quad \text{if}\quad |z_1| + |z_2| \geq 1. \endalign $$ 

Define $L_1 = L\circ \alpha.$ Then, for $D\subset \cc,$

$$\int_D L(\bd f, \be f) =  \int_D L_1 (\nabla f).$$

Thus, Conjecture 1 may be restated as: $L_1$ is quasiconvex at $0.$ 

In [Sv1], \v Sver\'ak introduced a class of functions containing $L_1$   
whose members he proved to be rank one convex, and noted that he was unable to 
determine if these functions are quasiconvex. We supply below a simple 
proof that $L_1$ is rank one convex.

For $A, B\in {\rr}^{2,2},$ write $\alpha(A) = (z_1,z_2),\;
\alpha (B) = (w_1, w_2).$ If $\text{rank}\, B\, \leq 1,$ then $|w_1| = |w_2|.$ 
Let $a=|z_1|^2 - |z_2|^2,\; b = 2\text{Re}\,(z_1 \overline{w_1}
- z_2 \overline{w_2}),$ and $I = \{t\in \rr: |z_1 + tw_1| + |z_2 + 
tw_2|<1\}.$ Then, for $g(t) = L_1(A+tB),$  we have,
when $\text{rank}\, B\leq 1,$
$$\align g(t) &= a+bt,\quad t\in I,\\
              &= 2|z_1 + tw_1| - 1, \quad t\in {\rr}\setminus I.\endalign$$

Now $g$ is continuous, $I$ is either empty or a bounded interval, and 
$t\rightarrow |z_1 + tw_1|$ is convex. It follows that $g$ is convex on 
$\rr.$ Hence, $L_1$ is rank one convex. Thus, if Conjecture 1 is false, then 
Morrey's conjecture for $m = n = 2$ will be confirmed.

It can also be shown that $L_1$ is not polyconvex. One way to do this is to show 
that $L_1$ does not satisfy condition (6) on [D1, p.107] when $A=0.$  

For $A\in {\rr}^{2,2},$ let $|A|^2 = a^2 + b^2 + c^2 + d^2.$ Let
$$E = \{A\in {\rr}^{2,2}: (|A|^2 +2\,\text{det\,A})^{1/2} 
+ (|A|^2 - 2\,\text{det\,A})^{1/2} \leq 2\}.$$

Then 
$$\align L_1(A) &= \text{det\,A},\quad A\in E, \\      
& = (|A|^2 +2\,\text{det\,A})^{1/2} - 1, \quad A\in {\rr}^{2,2}\setminus E.
 \endalign$$

Some rank one convex functions which look something like $L_1$ are studied in 
[DDGR] and [Sv3]. 

The connection between Morrey's conjecture and the Beurling-Ahlfors transform 
is discussed also in [As2] and [BL]. \bigskip 

\noindent 6. \bf Stretch Functions. \rm  There is a large class 
of functions for which equality holds in Conjectures 1 and 2. Write $z=r\eh.$ 
Functions $f:\cc \rightarrow \cc$ of the form  $$f(z) = g(r)\eh,$$
where $g$ is a nonnegative locally Lipschitz function on $(0,\infty)$ with 
$$g(0)  \equiv g(0+) = 0, \quad \text{and} \quad \lim_{r\rightarrow \infty} 
g(r) = 0,$$
will be called stretch functions. Let $\Cal S$ denote the set of all 
stretch functions. For $f\in \Cal S$, we have
$$\bd f = \hf(g' + r^{-1}g ),\quad \be f = \hf e^{2i\th}(g' - r^{-1}g),\quad
|\bd f| + |\be f| = \max(r^{-1}g, |g'|).\tag6.1$$

Let ${\Cal S}_1$ denote the subclass of $f\in \Cal S$ such that, for 
a.e. $r\in [0, \infty),$ holds $$|g'(r)|\leq r^{-1}g(r).\tag6.2$$

For example, for each $\alpha \in (0,1],\;\beta\in (0,1],$ and positive 
constant $c,$ the functions  
$$\align f(z) & = cr^{\alpha}\eh,\quad |z|\leq 1,\\
              & = cr^{-\beta}\eh,\quad |z|\geq 1, \tag6.3\endalign$$
belong to ${\Cal S}_1.$ 

\proclaim {Theorem 1} If $f\in {{\Cal S}_1} \cap \dot W^{1,2}(\cc, \cc),$ then 
$$\int_{\cc} L(\bd f, \be f) = 0.\tag6.4$$ 

If $f\in {{\Cal S}_1} \cap \dot 
W^{1,p}(\cc,\cc), \;$ then
$$\int_{\cc} {\Phi}_p(\bd f, \be f) = 0, \quad 1 < p \leq 2,\quad
  \int_{\cc} {\Phi}_p(\be f, \bd f) = 0,
\quad 2 \leq p < \infty.\tag6.5$$ \endproclaim

\demo {Proof} From (1.1a), (1.1b) and an approximation argument, it follows 
that we need only prove (6.4). Let $h(r) = g(r)/r.$ Then $h$ is  continuous 
and nonincreasing on $(0,\infty),$ with $\lim_{r\rightarrow \infty} h(r) = 0.$
From the last equation in (6.1) follows $|\bd f(z)| + |\be f(z)| = h(r).$ 

Let $E = \{r\in (0,\infty): h(r) > 1\}.$ Then 
$$\align L(\bd f(z), \be f(z)) & = r^{-1}g(r) + g'(r) - 1, \quad r\in E, \\
     & =r^{-1}g(r)g'(r), \quad  r\notin E.\tag6.6\endalign$$
  
Thus,
$$\align rL(\bd f(z), \be f(z)) & = {d\over {dr}}(rg - \hf r^2), \quad r\in E, 
\\
     & = \hf {d\over {dr}} g^2, \quad  r\notin E.\tag6.7\endalign$$

Now $E$ is either empty, or is a single interval $(0,R],$ with $0<R<\infty.$
Moreover, $g(0) = 0$ and $\lim_{r\rightarrow \infty} g(r) = 0.$ 
If $E$ is nonempty, then from (6.7) follows
$${1\over {2\pi}} \int_{\cc} L(\bd f, \be f) =   (Rg(R) -  \hf R^2) -  
\hf g^2(R) = -\hf (g(R) - R)^2.\tag6.8$$

The definition of $R$ implies that $g(R) = R.$ Hence (6.4) is true when $E$ is 
nonempty. If $E$ is empty, then it follows again from (6.7) that the integral 
on the left hand side of (6.8) equals $0.$ Theorem 1 is proved.  \enddemo

\proclaim {Theorem 2} If $f\in {\Cal S} \cap \dot W^{1,2}(\cc, \cc),$ then 
$$\int_{\cc} L(\bd f, \be f) \geq 0.\tag6.9$$ 

If $f\in {\Cal S} \cap \dot W^{1,p}(\cc,\cc), \;1<p<\infty,$ then
$$\int_{\cc} {\Phi}_p(\bd f, \be f) \geq 0.\tag6.10$$ \endproclaim

Thus, Conjectures 1 and 2 are true for stretch functions. According to Theorem 
1, the equality sign holds in (6.9) for the stretch functions which 
also satisfy (6.2). When $1<p\leq 2,$ the equality sign holds in (6.10) for 
stretch functions which satisfy  (6.2), while for $2\leq p<\infty$ equality
holds for their complex conjugates. 

As was pointed out to us by Iwaniec, by no means do all extremals for 
Conjectures 1 and 2 belong to ${\Cal S}_1.$  For example, start with the unit 
disk  
$B$ in the plane. For $j=1,2,...$, let $B_j = \{z: |z-a_j|<r_j\}$ 
be disjoint sub-disks of $B,$ and let $\{f_j\}$ be a sequence in 
${\Cal S}_1,$ with $f_j(z) = z$ on $|z| = 1.$ Define $f:\cc \rightarrow \cc$ by
$$\align f(z) & = a_j + r_j f_j({{z - a_j}\over {r_j}}), \quad \text{if}\quad 
z\in B_j,\\
 & = z,\quad \text{if} \quad z\in B\setminus \cup_{j=1}^{\infty} B_j, \\
 & = 1/\zb,\quad \text{if} \quad z\in \cc\setminus B.\endalign $$

Then, for each choice of $\{f_j\},$ equality holds for $f$ in Conjecture 1, and 
for $f$ in Conjecture 2 when $1<p\leq 2.$ Equality holds for $\overline{f}$ in 
Conjecture 2 when $2\leq p <\infty.$

Internal evidence, together with Burkholder's martingale results, suggests 
that when $p\ne 2$ there are no nontrivial functions for which equality is 
achieved in Conjecture 3. 

\demo {Proof of Theorem 2} As with Theorem 1, it suffices to prove (6.9). 
In our proof of (6.9) we shall assume that $g$ is continuously differentiable 
on $[0,\infty).$ The case when $g$ is locally Lipschitz then follows by an 
approximation argument. 

Let $E = \{r\in (0,\infty): |\bd f(z)| + |\be f(z)| > 1\}.$ Then, from (6.1),
$$\align L(\bd f(z), \be f(z)) & = |r^{-1}g(r) + g'(r)| - 1, \quad r\in E, \\
     & =r^{-1}g(r)g'(r), \quad  r\notin E.\tag6.11\endalign$$

For $r\in [0,\infty),$ define $F(r) = r^{-1}g + g' - 1.$ If $r\in E,$ then 
$$ F(r) \leq L(\bd f(z), \be f(z)),\tag6.12$$
by (6.11). If $r\notin E,$ then the third equation in (6.1) implies that
$g(r)\leq r$ and $|g'(r)|\leq 1.$ Hence, for $r\notin E,$
$$   F(r) - L(\bd f(z), \be f(z))  = r^{-1}g + g' - 1 - r^{-1}gg' 
=(1-r^{-1}g)(g'-1) \leq 0.$$

Thus, (6.12) holds for all $r\in [0,\infty).$  

If the set $\{r\in [0,\infty): g(r)\geq r\},$ is nonempty, let $R$ denote its 
supremum. If the set is empty, define $R=0.$ Then $0\leq R <\infty,$ since 
$g=o(1)$ at $\infty,$ and $g(R) = R.$ Since (6.12) holds for $r\in [0, 
\infty),$ we have 
$${1\over {2\pi}} \int_{|z|<R} L(\bd f, \be f) \geq  \int_0^R rF(r)\,dr =
(Rg(R) - \hf R^2) = \hf g^2(R).\tag6.13$$

Let $G(r) = L(\bd f(z), \be f(z)).$  Then $rG = gg'$ on $[0,\infty)\setminus 
E,$ by (6.11), and
$${1\over {2\pi}} \int_{|z|>R} L(\bd f, \be f)  
= \int_R^{\infty} gg'\,dr +  \int_{E\cap (R,\infty)} (rG - gg')\,dr.\tag6.14$$ 

From (6.13) and (6.14), it follows that
$${1\over {2\pi}} \int_{\cc} L(\bd f, \be f) \geq 
\int_{E\cap (R,\infty)} (rG - gg')\,dr.\tag6.15$$

Since $g(r) < r$ for $r>R,$ it follows from (6.1) that $E\cap (R,\infty) =
\{r\in (R,\infty): |g'(r)| > 1\}.$ If $E\cap (R,\infty)$ is empty, then (6.9) 
follows from (6.15). Assume $E\cap (R,\infty)$ is nonempty. Then it is a
finite or countable union of open intervals $(r_1, r_2)\subset (R, \infty),$
on each of which either $g'$ is everywhere $>1$ or $g'$ is everywhere $<-1.$
The hypothesis $f\in \dot W^{1,2}$ insures that all endpoints $r_2$ are finite. To 
prove (6.9), it suffices, in view of (6.15), to prove that, for each such
$(r_1, r_2),$
$$\int_{r_1}^{r_2} (rG-gg')\,dr \geq 0.\tag6.16$$

Suppose that $g'>1$ on $(r_1, r_2).$ Then, on $(r_1, r_2),$  
$$rG-gg' = g + rg' -r -gg' = - \hf {d\over {dr}}(r-g)^2.$$ 

Hence,
$$\int_{r_1}^{r_2} (rG-gg')\,dr = \hf [(r_1-g(r_1) )^2 - 
(r_2-g(r_2))^2].\tag6.17$$

But 
$$g'>1 \implies g(r_2) - g(r_1) > r_2 - r_1 \implies 
r_1 - g(r_1) > r_2 - g(r_2) > 0.$$

Thus, the integral in (6.17) is $>0.$

Suppose that $g'< -1$ on $(r_1, r_2).$ Then, on  $(r_1, r_2),$ 
$$rG-gg' = - g - rg' -r -gg' = - \hf {d\over {dr}}(r+g)^2.$$ 

Hence,
$$\int_{r_1}^{r_2} (rG-gg')\,dr = \hf [(r_1 + g(r_1) )^2 - 
(r_2 + g(r_2))^2].\tag6.18$$

But 
$$g'<-1 \implies g(r_2) - g(r_1) < - r_2 + r_1 \implies 
  r_1 + g(r_1) > r_2 + g(r_2) > 0.$$

Thus, the integral in (6.18) is $>0.$ The proof of (6.9) is complete.\enddemo
\bigskip

\noindent 7. \bf Some other partial results. \rm There are a few other 
classes of functions, in addition to the stretch functions, for which we 
can confirm Conjectures 1 and or 2. 

\proclaim {Theorem 3} For $a,b\in \cc,\; k = 1,2,3,...,$ Conjecture 1 is true 
for  
$$\align f(z) & = az^k + b{\overline{z}}^k, \quad |z| \leq 1, \\
   & = a{\overline{z}\,}^{-k} + bz^{-k}, \quad |z| \geq 1. \endalign$$ 
\endproclaim

\proclaim {Theorem 4} For $1<p<\infty,$ Conjecture 2 is true for $f\in 
\dot W^{1,p}(\cc, \cc)$  provided
$$f\quad \text{is harmonic in}\quad \cc \cup \{\infty\}\setminus \{|z| = 
1\},\tag7.1$$  
or
$$f=F\circ f_1\quad \text{or}\quad f = \overline{F} \circ f_1,\tag7.2$$ 
where $f_1\in \Cal S$ and $F$ is holomorphic on $f_1(\cc),$
\endproclaim  

Recall that $\Cal S$ denotes the class of stretch functions defined in Section 
6. The function $\Phi_p$ is homogeneous but the function $L$ is not; that is 
the main reason we can verify Conjecture 2 for more functions than 
we can for Conjecture 1. 

Theorem 3 can be proved by direct computation. The 
proof of Theorem 4 requires computation plus the fact that $p'th$ means 
of subharmonic functions on circles increase as the radius increases. We'll 
confine ourselves to sketching the proof of (7.1) when $p>2.$  

Suppose that $f\in \dot W^{1,p}(\cc)$ is harmonic in $|z|<1$ and in $1 < |z| \leq 
\infty.$ Then there 
exist holomorphic functions $g$ and $h$ in $|z|<1$ such that
$$\align f(z) &= g(z) + \overline{h}(z),\quad |z|<1,\\
          & = f(1/\overline{z}), \quad |z|>1.\endalign$$

Let $p>2$. Then computation gives
$$\align  {1\over {{\alpha}_p}}\int_{\cc} &{\Phi}_p(\bd f, \be f) 
 =  \int_{|z|<1}((p-1)|g'(z)| - |h'(z)|)(|g'(z)| + |h'(z)|)^{p-1}\,dx\,dy\\
& +  \int_{|z|>1} ((p-1)|h'(1/\overline{z})| - |g'(1/\overline{z})|)
(|h'(1/\overline{z})| + |g'(1/\overline{z})|)^{p-1} |z|^{-2p}\,dx\,dy\\
& = 2\pi \int_0^1((p-1) - r^{2p-4})I_1(r)\, r\,dr \,+\, 
2\pi \int_0^1 ((p-1)r^{2p-4} - 1) I_2(r)\,r\,dr,  \endalign$$
where $I_1(r),\;I_2(r)$ are the respective mean values on the circle $|z| = r$ 
of the functions $|g'| (|g'| + |h'|)^{p-1}$ and $|h'|(|g'| + |h'|)^{p-1}.$ The 
logarithms  
of these functions are subharmonic, hence so are the functions themselves. 
Thus, $I_1$ and $I_2$ are nondecreasing functions of $r$ on $[0,1].$
From $I_2 \nearrow,$ one easily shows that the integral containing $I_2$ is 
nonnegative. The integral containing $I_1$ is clearly nonnegative, because its 
integrand is. Hence, $\int_{\cc} {\Phi}_p(\bd f, \be f) \geq 0.$ \bigskip
                                               \bigskip
\def\T{\Bbb T}

\noindent 8. \bf Numerical Evidence. \rm
In this section, we present numerical evidence in favor of Conjecture~1.
Let $\T$ be the space $[0,1]$ with $0$ and $1$ identified.  Then
$W^{1,2}(\T^2,\cc)$ will denote the Sobolev space of complex valued functions
$f:[0,1]^2 \to \cc$ such that $f(0, y) \equiv f(1, y)$, $f(x, 0) \equiv f(x, 
1)$, and both $f$ and its distributional derivatives are in $L_2$. We will 
work with the following conjecture, which is equivalent to Conjecture~1.

\proclaim {Conjecture 4} Let $f\in W^{1,2}(\T^2, \cc).$ Then
$$\int_{\T^2} L(\bd f, \be f)
\geq 0. $$  \endproclaim

The approach is to consider piecewise linear functions described as
follows.  Let $N$ be a natural number.  Let $p_n$ be the fractional
part of $n/N$ (so that $p_{N+n} = p_n$).
Split $\T^2$ into triangles $\Delta^+_{m,n}$ with corners
$(p_m,p_n)$, $(p_{m+1},p_n)$, $(p_m,p_{n+1})$, and triangles
$\Delta^-_{m,n}$ with corners
$(p_m,p_n)$, $(p_{m-1},p_n)$, $(p_m,p_{n-1})$.

We will say that
$u:\T^2 \to \cc$ is an element of $\Cal P_N$ if $u$ is continuous, and
linear on
each of the triangles $\Delta^+_{m,n}$ and $\Delta^-_{m,n}$.  In this
way, once one knows that $u$ is an element of $\Cal P_N$, then $u$
is totally determined by its values at $(p_m,p_n)_{0\le m,n \le N-1}$.
Thus $\Cal P_N$ is a $2 N^2$ real dimensional space.  Let
$\iota : \rr^{2N^2} \to \Cal P_N$ denote an isomorphism.  Our goal
is to check whether the function $F_N:\rr^{2 N^2} \to \rr$ always takes
positive
values, where
$$ F_N(x) = \int_{\T^2} L(\bd (\iota x), \be (\iota x)) .$$
In fact, by an approximation argument, Conjecture~4 is equivalent to showing
that $F_N(x) \ge 0$ for all $x \in \rr^{2N^2}$ and all $N \ge 1$.

We obtained much numerical evidence to support this conjecture.
The algorithm was to choose a vector $x \in \rr^{2N^2}$ at random,
then minimize $F_N$, with $x$ as starting point, using the
conjugate gradient method described in Chapter~10.6 in [PTVF].  This was
done for various values of $N$, ranging from $6$ to $100$.  In every
case, it was found, up to machine precision, that $F_N$ always
takes non-negative values.  The results were verified independently
using Maple.

To implement this algorithm, it was necessary to compute
the gradient $\nabla F_N$. Because of the special nature of
this function, the computations needed to do this were not much more arduous
than the computations required for $F_N$.  The formulae required to
find $\nabla F_N$ were determined using Maple.

Other interesting facts emerged.  For a given $x \in \rr^{2N^2}$, we
may consider the function $h:\rr \to \rr$ given by
$$ h(t) = F_N(t x) .$$
It was found that this function is always increasing for $t \ge 0$, and
always decreasing for $t \le 0$.  However, it was also found
that the function $h$ is {\it not\/} necessarily convex.

This last fact is interesting, because if in Conjecture~4 the function
$L$ were to be replaced by a convex function, then $h$ would be convex.

 \enddocument